\documentclass{article} 

\usepackage{amsthm,amsmath,amssymb,amscd,enumerate,epsfig}
\usepackage{amsfonts}
\usepackage{epstopdf}
\usepackage[caption=false]{subfig}

\usepackage[numbers,sort&compress]{natbib}
\bibpunct[, ]{[}{]}{,}{n}{,}{,}
\renewcommand\bibfont{\fontsize{10}{12}\selectfont}

\theoremstyle{plain}
\newtheorem{theorem}{Theorem}[section]
\newtheorem{lemma}[theorem]{Lemma}

\newtheorem{observation}[theorem]{Observation}

\theoremstyle{definition}
\newtheorem{definition}[theorem]{Definition}
\newtheorem{example}[theorem]{Example}

\theoremstyle{remark}

\begin{document}

\title{Symbolic Dynamic Formulation for the Collatz Conjecture: I. Local and Quasi-global Behavior}

\author{
Eric Sakk \\ 
eric.sakk@morgan.edu \\
Department of Computer Science \\
Morgan State University, \\
Baltimore, USA.\\
}

\maketitle

\begin{abstract}
In this work, a symbolic dynamical formulation based upon discrete iterative mappings derived from the Collatz conjecture is introduced.
It is demonstrated that this formulation naturally induces a ternary alphabet useful for characterizing the expansive and dissipative behavior of generated itineraries.  Furthermore, local and quasi-global analyses indicate cyclic behaviors that should prove useful for
describing  global stability properties of itineraries. Additionally, techniques for generating arbitrarily long divergent itineraries are presented. Finally, this symbolic formulation allows itineraries to be grouped into sequence families that retain equivalent dynamical behavior.
\end{abstract}


 
\section{Introduction} \label{sec:intro}
In this work, we present a symbolic dynamical formulation for the \textbf{Collatz function}
\begin{equation} \label{eqn:colmain}
 f(x) =
 \left \{
 \begin{array}{ll}
 x/2 & x \; even \\
 3x+1 & x \; odd \\
 \end{array}
 \right.
\end{equation}
for any $x \in \mathbb{Z}^+$. Much has been written about this mapping \cite{lag1,lag2,lag3,lag4} and it has been observed (but, not rigorously proved) that all initial conditions tested so far converge to the periodic solution $1 \rightarrow 2 \rightarrow 4 \rightarrow 1 \rightarrow \cdots $ \cite{bar2021}. Intuitively, the above mapping should be describable using some type of shift map using integers represented in, for example, base 2 or base 3. Obviously, the $x/2$ portion of the mapping represents a right shift in base 2. The $3x+1$ portion represents a left shift plus one
in base 3; or, using $x + 2x +1$, the $3x+1$ portion could be described using a shift left plus $x$ plus $1$ in base 2. Since the Collatz function is a discrete dynamical system, we adopt terminology consistent with the analysis of this class of systems. For example. given an initial condition, the term '\textbf{\textit{itinerary}} ' is used to describe the set of points visited by iterating the Collatz function. 

The goal of this paper is to introduce
a symbolic formalism that enables the characterization of all possible
itineraries in order to determine if any given trajectory is unstable (i.e. divergent). Given that the Collatz function can be viewed as a piecewise series of left and right shifts, it seems sensible to suggest such an approach.  This work deviates from previous work concerning Collatz graphs \cite{Colussi2011,Hew2016} in that this formulation can be applied to characterize the long term behavior of families of trajectories. From the treatment presented here, it then becomes possible to identify dynamically equivalent symbolic sequences at many different scales. This work focuses on the local and quasi-global behavior of the proposed symbolic dynamical formulation. The results presented here are consistent with findings regarding the Collatz function \cite{terras1976,Jenber2021,schwob2021}. Similarities and differences will be pointed out in the upcoming sections.  

\section{Preliminary Concepts} \label{sec:prelims}
We begin by introducing  the definitions necessary for characterizing
the mapping defined in Equation (\ref{eqn:colmain}).

\subsection{Canonical Relation} \label{subsec:canonical}
A consequence of the Collatz function is that, given any positive, odd integer $x$, 
\begin{equation} \label{eqn:canonical}
f(x) = \frac{f(4x+1)}{4}.
\end{equation}
\noindent This relationship highlights a main feature of the formalism we introduce to characterize the dynamical behavior of the Collatz function. It also hints at the possibility of developing a renormalization procedure for understanding Equation (\ref{eqn:colmain}).

\subsection{The C-Ladder} \label{subsec:C-ladder}
Our goal then is to introduce an approach that can be used to order and collect the itineraries generated by any set of initial conditions. 
\begin{definition} \label{def:node}
A \textit{\textbf{node}} with respect to Equation (\ref{eqn:colmain}) is any $x \in \mathbb{Z}^+$ such that $x$ is odd.
\end{definition}

\begin{definition} \label{def:hub}
A \textit{\textbf{hub}} is a node that is not of the form $4k+1$ for any positive odd integer $k=1,3,5, \ldots$.
\end{definition}
\noindent In other words, all hubs are nodes, but not all nodes are hubs.

The hub is the starting point or 'generator' for an object we will refer to as a \textit{\textbf{C-ladder}}.  One side of the C-ladder is formed by starting with a hub and then \textit{\textbf{pumping the associated odd number}} $\bf{k}$ \textit{\textbf{using successive applications of}} $\bf{4k+1}$; hence, one side of the ladder consists of a sequence of nodes. In base 2, this sequence is formed by two left shifts then adding a value of one.  The 'rungs' of the ladder are generated by applying the Collatz function to the pumped nodes. The other side of a C-ladder is naturally induced by Equation (\ref{eqn:colmain}) and the canonical relation in Equation \ref{eqn:canonical}. For example, the odd integer $3$ is a hub and generates the following C-ladder 
\begin{equation*}
\begin{matrix}
3 & \Rightarrow & 13 & \Rightarrow & 53 & \Rightarrow & 213 & \Rightarrow & \cdots \\
\downarrow & & \downarrow & & \downarrow & & \downarrow & &  \\
10 & \leftarrow_4 & 40 & \leftarrow_4 & 160 & \leftarrow_4 & 640 & \leftarrow_4 & \cdots \\
\end{matrix}
\end{equation*}
In this diagram, numbers connected by $\rightarrow$ indicate a node and its image of the Collatz function. Numbers connected by $\leftarrow_4$ indicate division by 4. Numbers connected by $\Rightarrow$ indicate a $4k+1$ pump.  Given, for example, $x=53$, Equation \ref{eqn:canonical} tells us that $f(53) = f(4(53)+1)/4 = 160$.

Finally, we have the following
\begin{definition} \label{def:C-tree}
A \textit{\textbf{C-tree}} is the connected union of all C-ladders.
\end{definition}

Constructions similar to the even side of the C-ladder have been previously mentioned  \cite{Jenber2021,schwob2021}. This focal point makes sense as this side of the ladder 
is responsible for the convergence to a node under successive iterates of the Collatz function. 
When successive pumps are viewed in base 2, it is clear why this is so.  As successive $4k+1$ pumps of a node are applied to the C-ladder, the original node is shifted left by two bits and a one is then inserted into the least significant bit. After several pumps the suffix of the original value will look like $...01010101$. When the odd portion of the Collatz function is applied as $x+2x+1$, the pumped value $...01010101$ is added to the shifted value $...010101010$ leading to the values $....11111111$. To complete the $x+2x+1$ operation, simply add one to $....11111111$ and the resulting value looks like $.....00000000$ which will 'fall' down the even side of the C-ladder until the original node is reached. It is clear why some previous works have focused on this dissipative process in the hopes that some kind of convergence will ultimately result as unions of C-ladders are considered. While such an approach is worthy of exploration, there is no inherent guarantee that the final node value reached at the bottom of a particular C-ladder will be less than some previous iterate leading to one of the pumped nodes in the ladder.  If there were such a guarantee, the Collatz conjecture would have been solved.  On the other hand, C-ladders are extremely important as they represent strongly dissipative events along the itinerary of an initial condition.

\subsection{Node Definitions} \label{subsec:nodes}
We can say much more about the nodes of a C-ladder. In order to do so, the mapping in Equation (\ref{eqn:colmain}) will require us to partition the set of all nodes into three disjoint sets using the following definitions.
\begin{definition} \label{def:0node}
A \textit{\textbf{0-node}} is an odd integer of the form $6m+3$ for any nonnegative integer $m=0,1,2, \ldots$.
\end{definition}
\noindent Since $6m+3 = 3(2m+1)$, is should be clear that set of all 0-nodes represents the set of nonnegative odd integers that are divisible by three.
\begin{definition} \label{def:1node}
A \textit{\textbf{1-node}} is an odd integer of the form $6m+5$ for any nonnegative integer $m=0,1,2, \ldots$.
\end{definition}
\begin{definition} \label{def:2node}
A \textit{\textbf{2-node}} is an odd integer of the form $6m+1$ for any nonnegative integer $m=0,1,2, \ldots$.
\end{definition}
\noindent The set of 0-nodes, 1-nodes and 2-nodes comprise a disjoint partition of the set of nonnegative odd integers. \\

Finally, we require one more definition in order to further describe an intrinsic property of 0-nodes.
\begin{definition} \label{def:0nodechar}
Let $x$ be a $0-node$ of the form $6m+3$. The set of all 0-nodes can be further partitioned into a set defined by $3(6n+t)$ where $t=1,3 \; or \;5$
and $n=0,1,2,3, \ldots$. Given the partition $3(6n+t)$, the \textit{\textbf{character of a 0-node}} is the value $t \in \{1,3,5\}$ associated with any 0-node.
\end{definition}
\begin{example}
We can show all partitions of each character as follows: 
\begin{itemize}
\item 0-nodes of character $t=1$: $\{3, \text{\textbf{21}}, 39, 67, \ldots \}$
\item 0-nodes of character $t=3$: $\{9, 27, \text{\textbf{45}}, 63, \ldots \}$
\item 0-nodes of character $t=5$: $\{15, 33, 51, \text{\textbf{69}}, \ldots \}$
\end{itemize}
\end{example}
In this example, we have also highlighted nodes in boldface that cannot be hubs as they are of the form $4k+1$ where $k$ is odd. Finally, when describing a 0-node, the \textit{\textbf{character of the node is depicted with a superscript}} as follows: $3^{(1)}, \; 51^{(5)}, \; 27^{(3)}, \; etc$.

According to Equation (\ref{eqn:colmain}), backward iterations of 1-nodes and 2-nodes can lead to either even preimages or odd preimages. On the other hand, 0-nodes are unique as their backward iterations cannot possess odd preimages. Since a 0-node is of the form $3(2m+1)$, there cannot exist an odd integer $x$ such that $2^{-n}(3x+1)$ can map to $3y$ where $y$ is odd. Therefore, given a 0-node $x$, all Collatz preimages must be of the form $2^nx$. Thus, 0-nodes are the end point of successive divisions by 2. In some sense, the C-tree is 'surrounded' by 'threads' of even numbered itineraries of the Collatz function whose forward iterations ultimately converge to a 0-node on a Collatz ladder. To see this, consider the example from the previous section where the odd integer $3$ was a hub used to generate a C-ladder of the form
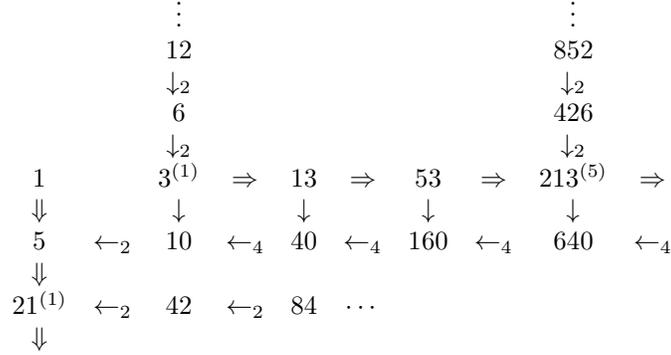
\begin{equation*}
\begin{matrix}
 &  & \vdots &  &   &   &   &   & \vdots  &   \\
 &  & 12 &   &   &   &   &  &   852 &   \\
 &  & \downarrow_2 &  &   &    &   &  & \downarrow_2  &   \\
 &  & 6 &   &   &   &   &   & 426  &   \\
 &  & \downarrow_2 &   &   &   &   &   & \downarrow_2   &   \\
1 &  & 3^{(1)} & \Rightarrow & 13 & \Rightarrow & 53 & \Rightarrow & 213^{(5)} & \Rightarrow \\
\Downarrow &  & \downarrow & & \downarrow & & \downarrow & & \downarrow & \\
5 &   \leftarrow_2 & 10 &  \leftarrow_4 & 40 & \leftarrow_4  & 160 & \leftarrow_4  & 640 & \leftarrow_4 \\
\Downarrow  &    &  &   &  &  &  &  &  &  \\
21^{(1)}  &   \leftarrow_2 & 42 &  \leftarrow_2 & 84 & \cdots  &   &   &  &   \\
\Downarrow  &   &   &   &  &  &   &  &  &  \\
\end{matrix}
\end{equation*}
(where $\leftarrow_4$ indicates division by 4 and $\leftarrow_2$ indicates division by 2).
This diagram also begins to hint at the fact that a C-ladder generated by a given hub can connect to another ladder generated by a different hub. Here, we can see that the 0-nodes 3, 21 and 213 are the end point of a series of successive divisions by 2. More importantly, we shall see in a moment that they exist at the beginning of what we will define as a '\textbf{\textit{primitive itinerary}}'. Before taking this step, we will use our node definitions to develop the symbolic formalism. Specifically, the node type will be included as part of the C-tree.
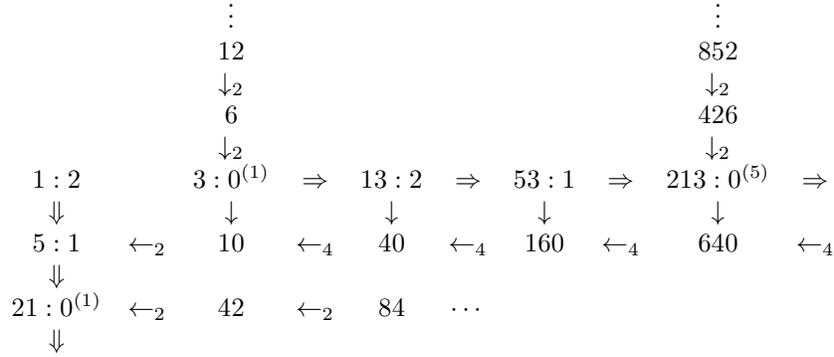
\begin{equation*}
\begin{matrix}
 &  & \vdots &  &   &   &   &   & \vdots  &   \\
 &  & 12 &   &   &   &   &  &   852 &   \\
 &  & \downarrow_2 &  &   &    &   &  & \downarrow_2  &   \\
 &  & 6 &   &   &   &   &   & 426  &   \\
 &  & \downarrow_2 &   &   &   &   &   & \downarrow_2   &   \\
1:2 &  & 3 : 0^{(1)} & \Rightarrow & 13:2 & \Rightarrow & 53:1 & \Rightarrow & 213: 0^{(5)} & \Rightarrow \\
\Downarrow &  & \downarrow & & \downarrow & & \downarrow & & \downarrow & \\
5:1 &   \leftarrow_2 & 10 &  \leftarrow_4 & 40 & \leftarrow_4  & 160 & \leftarrow_4  & 640 & \leftarrow_4 \\
\Downarrow  &    &  &   &  &  &  &  &  &  \\
21: 0^{(1)}  &   \leftarrow_2 & 42 &  \leftarrow_2 & 84 & \cdots  &   &   &  &   \\
\Downarrow  &   &   &   &  &  &   &  &  &  \\
\end{matrix}
\end{equation*}
According to this illustration, each node is directed to another node via successive application of the Collatz function. When the $4k+1$ pumps are viewed in base 2, the $f(x)=x+2x+1$ interpretation of the Collarz function makes clear how successive divisions by four down one side of a C-ladder can take place via a series of right shifts.  

\section{Local Behavior} \label{sec:local}
Given the above definitions, we now characterize the dynamical behavior of 1-nodes and 2-nodes. Consider the following examples.
\begin{equation*}
\begin{matrix}
& & 3:0 \\
  & & \downarrow \\
5:1 & \leftarrow_2 & 10 \\
\end{matrix}
 \quad  \quad  \quad \quad  \quad
\begin{matrix}
& & 17:1 \\
  & & \downarrow \\
13:2 & \leftarrow_4 & 52 \\
\end{matrix}
\end{equation*}
These cases indicate that a 1-node is the image of a hub under the mapping $(3x+1)/2$ and a 2-node is the image of a hub under the mapping $(3x+1)/4$. Even preimages can be recovered from a C-ladder using the canonical relation, so we will focus on hub images for a moment.

\begin{definition}
The \textit{\textbf{first odd preimage of a 1-node}} $y$ is formed from the inverse mapping $(2y-1)/3$.
\end{definition}

\begin{definition}
The \textit{\textbf{first odd preimage of a 2-node}} $y$ is formed from the inverse mapping $(4y-1)/3$.
\end{definition}

To establish the sensibility of these definitions requires the analyses presented in the following sections.

\subsection{Characterization of 1-nodes}
Given a node $x$, in order for $(3x+1)/2$ to map to a 1-node, the equation $(3(2m+1)+1)/2 = 6p+5$ must have a solution. This observation implies that there exists a solution to $6m+4 = 12p+10$ or $m = 2p+1$ for $p=0,1,2,\ldots$.  For example, $p=0$ leads to $6p+5 = 5$ and $m=2p+1 = 1 \Rightarrow 2m+1 = 3$; hence 3 maps to 5 which is a 1-node (see above diagram). It should also be clear that neither nodes of the form $6p+1$ (2-node) nor $6p+3$ (0-node) can be the image of $(3x+1)/2$. We have already dealt with the 0-node case. To understand the 2-node case, the equation $(3(2m+1)+1)/2 = 6p+1$  has no solution. Such a condition reduces to $6m+4 = 12p+2$ or $3m+1 = 6p$ which obviously cannot be satisfied for any integer $p$.

If we invert the 1-node value $6p+5$ to find its first odd preimage, we must solve $6p+5 = (3x+1)/2$; or $x = 4p+3$. Since $x$ cannot be of the form $4k+1$ with $k$ odd,  it must therefore, by definition, be a hub.
\begin{lemma}
The first odd preimage of a 1-node must always be a hub.
\end{lemma}
The hub itself can be a 0-node,, a 1-node or a 2-node, but the first odd preimage of a 1-node must be a hub.  Given the preceding analysis, we can also conclude that
\begin{lemma}
The first odd preimage of a 1-node of the form $6p+5$ must always be a hub of the form $4p+3$ where $p=0,1,2,3,\ldots$.
\end{lemma}

Given these observations, we have the following
\begin{theorem}
A hub $x$ of the form $4p+3$ under the image of the mapping $(3x+1)/2$ always maps to a 1-node.
\end{theorem}
Observe, the 1-node under this image of the hub need not also be a hub, but, as we shall see,  it can be.

We can collect these results in the form of a table to exemplify the behavior of 1-nodes.
\begin{table}[!ht]
\begin{center}
\caption{1-node behavior}
\label{table:1-node}
\begin{tabular}{|c|c|c|}
\hline
$p$ & (first odd preimage) $4p+3$: node type  & $6p+5$ 1-node \\
\hline
0 & $3:0^{(1)}$  & \textbf{5} \\
1 & $7:2$  & 11 \\
2 & $11:1$  & 17 \\
3 & $15:0^{(5)}$  & 23 \\
4 & $19:2$  & \textbf{29} \\
5 & $23:1$  & 35 \\
6 & $27:0^{(3)}$  & 41 \\
7 & $31:2$  & 47 \\
8 & $35:1$  & \textbf{53} \\
9 & $39:0^{(1)}$  & 59 \\
10 & $43:2$  & 65 \\
11 & $47:1$  & 71 \\
12 & $51:0^{(5)}$  & \textbf{77} \\
13 & $55:2$  & 83 \\
14 & $59:1$  & 89 \\
15 & $63:0^{(3)}$  & 95 \\
\hline
\end{tabular}
\end{center}
\end{table}
Table \ref{table:1-node} provides the first hint that, when viewed from the proper reference frame, there exist a number of interesting period behaviors in the set of generated 1-nodes. In addition to the preimage 0-node, 2-node, 1-node pattern as $p$ increases, for the set of 0-nodes, the type is periodic in the sequence $1,5,3$.  Furthermore, 1-nodes that are of the form $4k+1$ where $k$ is odd are indicated in boldface. Given their structure, there should be no surprise that these occur periodically as well.

The periodic nature of these values allows us to work backward from a 1-node value $y$ to determine node type $s$ of its first odd preimage using the formula
\begin{equation} \label{eqn:inverse_1-nodetype}
s = \frac{y+4}{3}  \quad (mod \; 3)
\end{equation}
which is the first hint at generating a symbolic dynamical description for the Collatz function.

\subsection{Characterization of 2-nodes}
Given a node $x$, in order for $(3x+1)/4$ to map to a 2-node, the equation $(3(2m+1)+1)/4 = 6p+1$ must have a solution. This observation implies that there exists a solution to $6m+4 = 24p+4$ or $m = 4p$ for $p=0,1,2,\ldots$.  For example, $p=0$ leads to $6p+1 = 1$ and $m=2p+1 = 1 \Rightarrow 2m+1 = 1$; hence 1 maps to 1 which is a 2-node. It should also be clear that neither nodes of the form $6p+5$ (1-node) nor $6p+3$ (0-node) can be the image of $(3x+1)/4$. We have already dealt with the 0-node case. To understand the 1-node case, the equation $(3(2m+1)+1)/4 = 6p+5$  has no solution. Such a condition reduces to $6m+4 = 24p+20$ or $6m= 24p+16$ which cannot be satisfied for any integer $p$.

If we invert the 2-node value $6p+1$ to find its first odd preimage, we must solve $6p+1 = (3x+1)/4$; or $x = 8p+1$. Since $x$ cannot be of the form $4k+1$ with $k$ odd,  it must therefore, by definition, be a hub.
\begin{lemma}
The first odd preimage of a 2-node must always be a hub.
\end{lemma}
The hub itself can be a 0-node,, a 1-node or a 2-node, but the first odd preimage of a 2-node must itself be a hub.  Given the preceding analysis, we can also conclude that
\begin{lemma}
The first odd preimage of a 2-node of the form $6p+1$ must always be a hub of the form $8p+1$ where $p=0,1,2,3,\ldots$.
\end{lemma}

Given these observations, we have the following
\begin{theorem}
A hub $x$ of the form $8p+1$ under the image of the mapping $(3x+1)/4$ always maps to a 2-node.
\end{theorem}
Observe, the 2-node under this image of the hub need not also be a hub, but, as we shall see,  it can be.

We can collect these results in the form of a table to exemplify the behavior of 2-nodes.
\begin{table}[!h]
\begin{center}
\caption{2-node behavior}
\label{table:2-node}
\begin{tabular}{|c|c|c|}
\hline
$p$ & (first odd preimage) $8p+1$: node type  & $6p+1$ 2-node \\
\hline
0 & $1:2$  & 1 \\
1 & $9:0^{(3)}$  & 7 \\
2 & $17:1$  & \textbf{13} \\
3 & $25:2$  & 19 \\
4 & $33:0^{(5)}$  & 25 \\
5 & $41:1$  & 31 \\
6 & $49:2$  & \textbf{37} \\
7 & $57:0^{(1)}$  & 43 \\
8 & $65:1$  & 49 \\
9 & $73:2$  & 55 \\
10 & $81:0^{(3)}$  & \textbf{61} \\
11 & $89:1$  & 67 \\
12 & $97:2$  & 73 \\
13 & $105:0^{(5)}$  & 79 \\
14 & $113:1$  & \textbf{85} \\
15 & $121:2$  & 91 \\
\hline
\end{tabular}
\end{center}
\end{table}
Table \ref{table:2-node} provides the first hint that, when viewed from the proper reference frame, there exist a number of interesting period behaviors in the set of generated 2-nodes. In addition to the preimage 2-node, 0-node, 1-node pattern as $p$ increases, for the set of 0-nodes, the type is periodic in the sequence $3,5,1$.  Furthermore, 2-nodes that are of the form $4k+1$ where $k$ is odd are indicated in boldface. Given their structure, there should be no surprise that these occur periodically as well.

The periodic nature of these values allows us to work backward from a 2-node value $y$ to determine node type $s$ of its first odd preimage using the formula
\begin{equation} \label{eqn:inverse_2-nodetype}
s = \frac{2y+4}{3}    \quad (mod \; 3)
\end{equation}
which is the another hint at generating a symbolic dynamical description for the Collatz function

\subsection{Characterization of the set of all nodes}
As pointed out above, even preimages of a hub image can be recovered from a C-ladder using the canonical relation by applying successive pumps of the mapping $4k+1$ to the hub.  Therefore, a C-ladder can be generated using 0-node, 1-node, 2-node as a hub and then pumping that value using $4k+1$.  From the previous sections, all hubs are either of the form $4p+3$ (which maps to a 1-node) or $8p+1$ (which maps to a 2-node).  Given these possibilities, if we partition of all odd numbers using the rule $8p+t$ where $t \in \{1,3,5,7\}$, we can see that the set of all odd numbers consists of only hubs and 'pumps'. The value $8p+1$ is a hub whose image corresponds to a 2-node. If we define, $p=2q$ is even and $p+1 = 2q+1$ is odd, then $4p+3$ (a hub that maps to a 1-node) corresponds to $8q+3$ and $4(2q+1)+3$ corresponds to $8p+7$; while $4(2p+1)+1$ pumps correspond to $8p+5$. Therefore, using the partition of all odd numbers using the rule $8p+t$ where $t \in \{1,3,5,7\}$ leads to 1-node hubs when $t \in \{3,7\}$, 2-node hubs when $t=1$ and pumps when $t=5$.

This characterization is significant because it tells us the how dissipative and expansive nodes are distributed. A hub that transitions to a 1-node is expansive as it is of the form $(3x+1)/2$. A hub that transitions to a 2-node is dissipative as it is of the form $(3x+1)/4$. A pumped $4k+1$ node (where $k$ is odd) is dissipative by at least a factor of $(3x+1)/8$ if it transitions to a 1-node and at least a factor of $(3x+1)/16$ if it transitions to a 2-node. These observations will allow us to begin entering into a discussion regarding the stability of trajectories.

\section{Quasi-global Behavior} \label{sec:quasi-global}

\subsection{Symbolic construction and the stability of an itinerary} \label{subsec:symostability}
In order to begin a discussion about stability, we require the definition
\begin{definition}
A \textit{\textbf{primitive itinerary}}  of Equation (\ref{eqn:colmain}) is one that begins at a 0-node and ends at a pumped $4k+1$ node where $k$ is odd.
\end{definition}
\begin{example}
$21:0$ is a primitive itinerary of length one as it is both a 0-node and a pumped $4k+1$ node where $k=5$.
\end{example}
\begin{example}
$3:0 \rightarrow_2 5:1$ is a primitive itinerary of length two.
\end{example}
\begin{example}
$81:0 \rightarrow_4 61:2 $ is a primitive itinerary of length two.
\end{example}
\begin{example}
$15:0 \rightarrow_2 23:1  \rightarrow_2 35:1 \rightarrow_2 53:1 $ is a primitive itinerary of length four.
\end{example}

More directly, a primitive itinerary begins at a 0-node and \textit{\textbf{traverses a sequence of hubs}} that are either 1-nodes or 2-nodes until a pumped node is encountered. A primitive itinerary will \textit{\textbf{naturally induce a sequence of 1's and 2's}}.
\begin{example}
$15:0 \rightarrow_2 23:1  \rightarrow_2 35:1 \rightarrow_2 53:1 $ is a primitive itinerary that generates the sequence    $ \{0111\}$.
\end{example}
\begin{example}
$9:0 \rightarrow_4 7:2  \rightarrow_2 11:1 \rightarrow_2 17:1 \rightarrow_4 13:2 $ is a primitive itinerary
that generates the sequence $ \{02112\}$.
\end{example}

If the value of the pumped node is less than that of the 0-node, then we say the primitive itinerary is \textit{\textbf{dissipative}}. If the value of the pumped node is greater than that of the 0-node, then we say the primitive itinerary is \textit{\textbf{expansive}}. We can make an asymptotic statement regarding the dissipative or expansive behavior of a primitive itinerary. After the 0-node, a primitive itinerary will consist only of 1-nodes and 2-nodes (i.e a sequence of 1's and 2's). It is therefore possible to asymptotically characterize the stability of such an itinerary by counting the occurrences of 1's and 2's in a sequence generated by the itinerary.

\begin{observation}
Letting $n$ represent the number of 2's and $m$ represent the number of 1's, a primitive itinerary will be expansive if the following condition holds:
\begin{equation}
\left(\frac{3}{2}\right)^m\left(\frac{3}{4}\right)^n > 1  \Rightarrow m >n \frac{\log{4/3}}{\log{3/2}} \approx 0.7095n \equiv \gamma n
\end{equation}
\end{observation}
This bound is consistent results regarding the Collatz function (see for example, \cite{lag5}).

After the 0-node, a primitive itinerary is simply a binary sequence of length $N$ and the binomial theorem can be used to count the number of expansive sequences $N_E$ generated by a primitive itinerary:
\begin{equation}
N_E = \sum_{k=0}^{\hat{n}} {N \choose k}
\end{equation}
where $\hat{m} > \gamma \hat{n}$ and $\hat{m} = N - \hat{n}$. This, in turn, yields the condition that
\begin{equation}
N - \hat{n} > \gamma \hat{n}  \Rightarrow \hat{n} < \frac{N}{1+ \gamma}
\end{equation}
meaning that\textit{\textbf{ the number of 2's in an expansive primitive sequence (i.e. a sequence generated by a primitive itinerary) of length $\mathbf{N}$ asymptotically must be less than the quantity }} $\mathbf{N/(1 + \gamma)}$.
\begin{example}
Consider primitive sequences of length $N=5$ with the first symbol equal $0$. Then, $\hat{n} < \frac{N}{1+ \gamma} \approx 2.3$. Therefore, expansive sequences can visit at most $2$ 2-nodes.
\end{example}

A complete itinerary of an initial condition is simply the intersection of a set of primitive itineraries. Once a primitive itinerary lands at a $4k+1$ node, it will necessarily drop by at least a factor of $8$ into the midst of another primitive itinerary and so on. This behavior is readily observable if one plots each value as function of its iteration index. For complex trajectories, it is common to see an expansive portion that suddenly dissipates and intersects another primitive itinerary and then recommences the expansion-dissipation process. Obviously, in order to encompass the possibility of returning to the $1 \rightarrow 2 \rightarrow 4$ cycle, there must ultimately be more dissipative events than expansive events.

\subsubsection{An expansive construction} \label{subsec:expansive}
One major goal of Section \ref{sec:quasi-global} will be to systematically characterize exactly how to generate expansive sequences as this would deal directly with the question regarding the existence of an asymptotically expansive sequence. Before presenting a more comprehensive development, we can deal more directly with the worst case scenario, \textit{\textbf{a sequence consisting purely of 1-nodes that are hubs}}.

Table \ref{table:1-node} indicates that it should be possible generate a sequence consisting purely of hubs that are of 1-nodes since 1-nodes that map to 1-nodes are separated by a value of 18. Therefore, given a 1-node $y$, there should exist an $n$ such that $y+18n = (3y+1)/2$. We can begin exploring this goal by checking that a hub can be mapped to a 1-node using the condition
\begin{equation} \label{eqn:1-nodeto1-node}
[3(4p+3)+1]/2 = 6p+5
\end{equation}
To ensure this condition can be iterated to generate a sequence of hubs that are 1-nodes, the following observation will prove useful.
\begin{observation} \label{obs:1-nodeiterator}
Given a node $y_i$ and the mapping $y_{i+1}=(3y_i+1)/2$
\begin{equation}
y_{i+1} + 1 =\frac{3y_i+1}{2} + 1 = \frac{3}{2}(y_i+1)
\end{equation}
\end{observation}
The procedure is to start with an $x$ that maps to 1-node $y_0$ that is a hub. Given Equation (\ref{eqn:1-nodeto1-node}), it follows that
\begin{equation} \label{eqn:1-nodeto1-node_ic}
y_0 + 1 = 6p+6 = 6(p+1).
\end{equation}
Observation \ref{obs:1-nodeiterator} tells us that 1-node iterations can be constructed simply by successively multiplying $y_i+1$ by an expansive factor of $3/2$.  In order to ensure that $y_i+1$ continues to be an integer after successive divisions by $2$, we must choose $p=s2^M -1$, for some $s$ an odd integer, in Equation (\ref{eqn:1-nodeto1-node_ic}) resulting in
\begin{equation} \label{eqn:1-nodeto1-node_ic2}
y_0 + 1 = 6s2^M = 3s2^{M+1}.
\end{equation}
For this prescription, it should be clear that $p$ must be an odd integer in order to create the initial condition. It should also be clear that $\bf{M}$ c\textit{\textbf{an be chosen arbitrarily large in order to create an arbitrarily long subsequence of}} $\bf{M+1}$ \textit{\textbf{1-nodes}} within a primitive sequence by applying Observation \ref{obs:1-nodeiterator} to Equation (\ref{eqn:1-nodeto1-node_ic2}) until all powers of two have been exhausted
\begin{equation}
\begin{split}
y_0 + 1 & = 3s2^{M+1} \\
y_1 + 1 & = \left( \frac{3}{2}\right) (3s2^{M+1}) = s3^2 2^M \\
y_2 + 1 & = \left( \frac{3}{2}\right) (s3^2 2^M) = s3^3 2^{M-1} \\
 & \vdots \\
y_{M+1} + 1 & = \left( \frac{3}{2}\right) (s3^{M+1} 2) =  s3^{M+2}\\
\end{split}
\end{equation}
The value $y_{M+1}$ will necessarily be an even integer indicating the last iteration can map to a 2-node or a pumped node.
\begin{example}
Consider the case where $s=1$, $M=5$ and $p=2^5-1 = 31$; hence,
\begin{equation*}
\begin{split}
y_0 + 1 & = 3s2^{M+1} = 192 \\
y_1 + 1 & = \left( \frac{3}{2}\right) (192) = 288 \\
y_2 + 1 & = \left( \frac{3}{2}\right) (288) = 432 \\
y_3 + 1 & = \left( \frac{3}{2}\right) (432) = 648 \\
y_4 + 1 & = \left( \frac{3}{2}\right) (648) = 972 \\
y_5 + 1 & = \left( \frac{3}{2}\right) (972) = 1458 \\
y_6 + 1 & = \left( \frac{3}{2}\right) (1458) = 2187 \\
\end{split}
\end{equation*}
This implies the following sequence containing $M+1=6$ successive 1-nodes that are hubs (and, hence, a subsequence within a primitive sequence),
\begin{equation*}
191:1 \rightarrow_2 287:1 \rightarrow_2 287:1 \rightarrow_2 431:1 \rightarrow_2 647:1 \rightarrow_2 971:1 \rightarrow_2 1457:1 \rightarrow_2 2188
\end{equation*}
The end point of this sequence is an even integer that can still be further divided by 2 (this is because the end of the sequence $1457$ maps to a 2-node).
\end{example}

To check that the prescribed 1-node to 1-node mapping truly results in a new 1-node we first check the following. Given a value $6p+5$ for some $p$,
there must exist a $q$ such that $6p+5 = 4q+3$, it is easily found that, when $p$ is odd, there must exist a $q$. But, recall, the condition $p =s2^M-1$ will be odd if $s$ is odd. Therefore, the initial condition must represent a hub as it maps a 1-node to a 1-node. Next, given Observation \ref{obs:1-nodeiterator} to Equation (\ref{eqn:1-nodeto1-node_ic2}), it is enough to check the first iteration on the initial condition
\begin{equation}
\begin{split}
y_1 + 1 & = (\frac{3}{2}) (y_0 + 1) \\
& = (\frac{3}{2}) (6p+6) \\
& = 9p+9 \\
& = 9(s2^M-1)+ 9 \\
\Rightarrow y_1 & = 9s2^M -1  \\
\end{split}
\end{equation}
Rearranging and regrouping $y_1$ as follows yields the desired result
\begin{equation}
\begin{split}
y_1  & = 9s2^M -1 \\
& = 3 \cdot 2 s2^{M-1} -6 + 5 \\
& = 6(s2^{M-1}-1) + 5 \\
& = 6q+5 \\
\end{split}
\end{equation}
which is of the form of a 1-node.
Hence, the initial prescription $p =s2^M-1$ does enable the mapping of 1-nodes to 1-nodes a total of $M$ times.

Given the possibility of generating an expansive subsequence of 1-nodes of arbitrary length, it is only natural to attempt to enumerate and categorize the set of all possible expansive primitive itineraries. This is the overall goal of the following subsections.

\subsection{Cyclic Behavior}
Given the characterization of nodes into the partition $\{6p+1, 6p+3, 6p+5\}$, it should not be surprising that the C-ladder exhibits a substantial amount of cyclic behavior. We state many of these observations without proof as they are fairly straightforward to demonstrate.
\begin{observation} \label{obs:pumpper1}
The node sequence of $4k+1$ pumps along a C-ladder obey the periodic node ordering
\begin{equation}
\cdots \Rightarrow 0 \Rightarrow 2 \Rightarrow 1 \Rightarrow 0 \Rightarrow 2 \Rightarrow 1 \Rightarrow \cdots.
\end{equation}
\end{observation}
We can say a bit more, given the character of a 0-node.
\begin{observation} \label{obs:pumpper2}
The node sequence of $4k+1$ pumps along a C-ladder obey the periodic node ordering
\begin{equation}
\begin{split}
\cdots\Rightarrow &  0^{(1)}  \Rightarrow 2 \Rightarrow 1 \Rightarrow \\
& 0^{(5)} \Rightarrow 2 \Rightarrow 1 \Rightarrow \\
& 0^{(3)} \Rightarrow 2 \Rightarrow 1 \Rightarrow \\
&  0^{(1)}  \Rightarrow 2 \Rightarrow 1 \Rightarrow \\
& 0^{(5)} \Rightarrow 2 \Rightarrow 1 \Rightarrow \\
& 0^{(3)} \Rightarrow 2 \Rightarrow 1 \Rightarrow \cdots \\
\end{split}
\end{equation}
where the character of 0-nodes is also periodic.
\end{observation}
As 0-nodes can possess different characters, the above observation begins to hint at the possibility of topologically equivalent behaviors between different C-ladders in a C-tree. To explore this possibility further, it is important to recognize that C-ladders are numerically constrained along its sequence of $4k+1$ pumps,
\begin{observation} \label{obs:pumpper3}
Consider a node sequence of $4k+1$ pumps along a C-ladder. Along the pump, the first odd preimage of 1-nodes and 2-nodes obey the periodic node ordering
\begin{equation}
\begin{split}
\cdots \Rightarrow &  0^{(1)}  \Rightarrow (1 \rightarrow_4 2) \Rightarrow (1 \rightarrow_2 1) \Rightarrow \\
& 0^{(5)} \Rightarrow (0 \rightarrow_4 2) \Rightarrow (2 \rightarrow_2 1) \Rightarrow \\
& 0^{(3)} \Rightarrow (2 \rightarrow_4 2) \Rightarrow (0 \rightarrow_2 1) \Rightarrow \\
& 0^{(1)}  \Rightarrow (1 \rightarrow_4 2) \Rightarrow (1 \rightarrow_2 1) \Rightarrow \\
& 0^{(5)} \Rightarrow (0 \rightarrow_4 2) \Rightarrow (2 \rightarrow_2 1) \Rightarrow \\
& 0^{(3)} \Rightarrow (2 \rightarrow_4 2) \Rightarrow (0 \rightarrow_2 1) \Rightarrow  \cdots \\
\end{split}
\end{equation}
\end{observation}

\begin{table}[!h]
\begin{center}
\caption{C-ladder successive $4k+1$ pump behavior starting with the hub value of 1 (*the hub with periodic itinerary $1 \rightarrow 2 \rightarrow 4$ is not primitive).}
\label{table:k-pumps1}
\begin{tabular}{|c|c|l|c|l|}
\hline
0-node value & 2-node value & primitive sequence & 1-node value & primitive sequence \\
\hline
-- & $1^\ast$  &  $222 \ldots$ & $5$  &  $10^{3(1)}$ \\
$21^{(1)}$ & $85$  &  $210^{75(1)}$ & $341$  &  $1120^{201(1)}$ \\
$1365^{(5)}$ & $5461$  &  $20^{7281(3)}$ & $21845$  &  $12210^{17259(3)}$ \\
$87381^{(3)}$ & $349525$  &  $2211120^{24581(5)}$ & $1398101$  &  $10^{932067(3)}$ \\
$5592405^{(1)}$ & $22369621$  &  $21120^{17674761(1)}$ & $89478485$  &  $1110^{26512143(5)}$ \\
$357913941^{(5)}$ & $1431655765$  &  $20^{(1)}$ & $5726623061$  &  $120^{(3)}$ \\
\hline
\end{tabular}
\end{center}
\end{table}

We can collect the above observations in the form of a table. In Table \ref{table:k-pumps1}, in the first column, the character of each 0-node value is included. The third column provides the complete sequence leading to the 1-node value in column two, all sequences therefore begin at a 0-node in the rightmost position, travel to the left and end at the 1-node value. The fifth column provides the complete sequence leading to the 2-node value in column four, all sequences therefore begin at a 0-node in the rightmost position, travel to the left and end at the 2-node value. In the sequence columns, for convenience, the character and value of the starting 0-node is provided. Because the 1-node and 2-node columns represent $4k+1$ pumps and all sequences begin at 0-nodes, by definition, all sequences (except the hub with periodic itinerary $1 \rightarrow 2 \rightarrow 4$ ) represent primitive itineraries. Table \ref{table:k-pumps3} shows a similar sequence analysis using $3$ as the hub.

\begin{table}[!h]
\begin{center}
\caption{C-ladder successive $4k+1$ pump behavior starting with the hub value of 3.}
\label{table:k-pumps3}
\begin{tabular}{|c|c|l|c|l|}
\hline
0-node value & 2-node value & primitive sequence & 1-node value & primitive sequence \\
\hline
$3^{(1)}$ & $13$  &  $21120^{9(3)}$ & $53$  &  $1110^{15(5)}$ \\
$213^{(5)}$ & $853$  &  $\bf{20^{1137(1)}}$ & $3413$  &  $\bf{120^{3033(3)}}$ \\
$13653^{(3)}$ & $54613$  &  $\bf{220^{97089(5)}}$ & $218453$  &  $\bf{10^{145635(5)}}$ \\
$873813^{(1)}$ & $3495253$  &  $\bf{212220^{(5)}}$ & $13981013$  &  $\bf{110^{(1)}}$ \\
$55924053^{(5)}$ & $223696213$  &  $20^{(5)}$ & $894784853$  &  $1210^{(3)}$ \\
$3579139413^{(3)}$ & $14316557653$  &  $222120^{(3)}$ & $57266230613$  &  $10^{(1)}$ \\
$0^{(1)}$ &    &  $210^{(3)}$ &    &  $112210^{(3)}$ \\
$0^{(5)}$ &    &  $20^{(3)}$ &    &  $12220^{(5)}$ \\
\hline
\end{tabular}
\end{center}
\end{table}

In addition to confirming Observations \ref{obs:pumpper1}-\ref{obs:pumpper3}, several other observations are possible. For example, as one steps down the rows in the sequence column, the symbolic structure induces a form of counting in each sequence position over the alphabet $\{0,1,2\}$. There is also much more to say regarding character periodicities. In fact, given this symbolic structure, C-trees and C-ladders exhibit a wide spectrum of periodic behaviors and topological consistencies that extend globally throughout the tree. Let us take another look at these behaviors from the perspective of the character of 0-nodes.
\begin{table}
\begin{center}
\caption{In each row of this table, each 0-node value generated by $6m+3$ (and starting with $m=0$) is followed by its 2-node and 1-node generated by $4k+1$ pumps.}
\label{table:0-nodetopology}
\begin{tabular}{|c|c|l|c|l|}
\hline
0-node & 2-node  & sequence & 1-node  & sequence \\
\hline
$3^{(1)}$ & $13$  &  $21120^{9(3)}$ & $53$  &  $1110^{15(5)}$ \\
$9^{(3)}$ & $37$  &  $22120^{57(1)}$ & $149$  &  $10^{99(3)}$ \\
$15^{(5)}$ & $61$  &  $20^{81(3)}$ & $245$  &  $122220^{513(3)}$ \\
$\bf{21^{(1)}}$ & $85$  &  $210^{75(1)}$ & $341$  &  $1120^{201(1)}$ \\
$27^{(3)}$ & $109$  &  $22210^{171(3)}$ & $437$  &  $10^{291(1)}$ \\
$33^{(5)}$ & $133$  &  $20^{177(5)}$ & $533$  &  $1210^{315(3)}$ \\
$39^{(1)}$ & $157$  &  $21210^{123(5)}$ & $629$  &  $110^{279(3)}$ \\
$\bf{45^{(3)}}$ & $181$  &  $220^{321(5)}$ & $725$  &  $10^{483(5)}$ \\
$\bf{51^{(5)}}$ & $\bf{205}$  &  $\bf{20^{273(1)}}$ & $\bf{821}$  &  $\bf{120^{729(3)}}$ \\
$57^{(1)}$ & $229$  &  $2110^{135(3)}$ & $917$  &  $111222121212220^{897(5)}$ \\
$63^{(3)}$ & $253$  &  $2211221220^{417(1)}$ & $1013$  &  $10^{675(3)}$ \\
$\bf{69^{(5)}}$ & $277$  &  $20^{369(3)}$ & $1109$  &  $1221120^{777(1)}$ \\
$75^{(1)}$ & $301$  &  $210^{267(5)}$ & $1205$  &  $112120^{633(1)}$ \\
$81^{(3)}$ & $325$  &  $22221110^{303(5)}$ & $1301$  &  $10^{867(1)}$ \\
$87^{(5)}$ & $348$  &  $20^{465(5)}$ & $1397$  &  $1211120^{489(1)}$ \\
$\bf{93^{(1)}}$ & $373$  &  $2120^{441(3)}$ & $1493$  &  $110^{663(5)}$ \\
$99^{(3)}$ & $397$  &  $220^{705(1)}$ & $1589$  &  $10^{1059(5)}$ \\
$105^{(5)}$ & $421$  &  $20^{561(3)}$ & $1685$  &  $120^{1497(1)}$ \\
$111^{(1)}$ & $445$  &  $211121121221111210^{27(3)}$ & $1781$  &  $11110^{351(3)}$ \\
$\bf{117^{(3)}}$ & $469$  &  $2210^{555(5)}$ & $1877$  &  $10^{1251(3)}$ \\
$123^{(5)}$ & $493$  &  $20^{657(3)}$ & $1973$  &  $1220^{2357(5)}$ \\
$129^{(1)}$ & $517$  &  $210^{459(3)}$ & $2069$  &  $1122211212120^{1017(3)}$ \\
\hline
\end{tabular}
\end{center}
\end{table}

Once again, in Table \ref{table:0-nodetopology}, in addition to periodic behavior, a ordered 'counting' process takes place in the sequences associated with each primitive itinerary (since these sequences are constrained by Observations \ref{obs:pumpper1}-Observations \ref{obs:pumpper1}). Given the $\{6m+3\}$ ordering of the 0-nodes, the second position of each sequence generates the periodic ordering $\{0,1,2\}$ for the 2-node column and $\{0,2,1\}$ for the 1-node column for every change by $6$ in the 0-node value. The third position of each sequence generates the periodic ordering $\{0,2,1\}$ for the 2-node column and $\{0,1,2\}$ for the 1-node column for every change by $18$ in the 0-node value. The fourth position of each sequence generates the periodic ordering $\{0,1,2\}$ for the 2-node column and $\{0,2,1\}$ for the 1-node column for every change by $54$ in the 0-node value and so on.

These observations imply an immense amount of consistent topological structure that is globally applicable across the C-tree. To see this, it is best to start with a motif that generates short primitive sequences that terminate quickly for most sequences in the motif. For example, consider the block of sequences highlighted in Table \ref{table:k-pumps3} that pump across three iterations of $4k+1$. The start of this motif can be found in Table \ref{table:0-nodetopology} at the 0-node labeled $51$. By simply adding 54 consecutively ($51, 105, 159, \ldots$), the reader can check that a measure of topological consistency and counting periodicity is preserved.  To preserve even greater topological consistency, one can take into account and preserve the character of the 0-node generating the primitive sequence (as these exhibit periodic behavior as well). More will be said about these behaviors when we present the global treatment in part II of this work. Before tackling the global behavior, we require one more set of observations.

\begin{table}
\begin{center}
\caption{Characterization of prefixes within a sequence ending at a node. Here, it is convenient to use the 1-node and 2-node definitions
given in Section \ref{subsec:nodes} in order to generate a natural ordering. Boldface sequences are primitive sequences.}
\label{table:prefixes}
\begin{tabular}{|c|c|l|c|l|}
\hline
$p$ & 2-node 6p+1 & sequence & 1-node 6p+5 & sequence \\
\hline
$0$ & $1_{\leftarrow_4} 1$  &  $222 \ldots$
        & $\bf{5_{\leftarrow_2} 3}$  &  $\bf{10^{(1)}}$ \\
$1$ & $7_{\leftarrow_4} 9$   &  $20^{(3)}$
        & $11_{\leftarrow_2} 7_{\leftarrow_4} 9$  &  $120^{(3)}$ \\
$2$ & $\bf{13_{\leftarrow_4}  17_{\leftarrow_2} 11_{\leftarrow_2} 7_{\leftarrow_4} 9}$  &  $\bf{21120^{(3)}}$
        & $17_{\leftarrow_2} 11_{\leftarrow_2} 7_{\leftarrow_4} 9$  &  $1120^{(3)}$ \\
$3$ & $19_{\leftarrow_4}  25_{\leftarrow_4} 33$  &  $220^{(5)}$
        & $23_{\leftarrow_2} 15$  &  $10^{(5)}$ \\
$4$ & $25_{\leftarrow_4} 33$  &  $20^{(5)}$
        & $\bf{29_{\leftarrow_2} 19_{\leftarrow_4} 25_{\leftarrow_4} 33}$  &  $\bf{1220^{(5)}}$ \\
$5$ & $31_{\leftarrow_4} 41_{\leftarrow_2} 27$  &  $210^{(3)}$
        & $35_{\leftarrow_2} 23_{\leftarrow_2} 15$  &  $110^{(5)}$ \\
$6$ & $\bf{37_{\leftarrow_4}  49_{\leftarrow_4} 65_{\leftarrow_2} 43_{\leftarrow_4} 57}$  &  $\bf{22120^{(1)}}$
        & $41_{\leftarrow_2} 27$  &  $10^{(3)}$ \\
$7$ & $43_{\leftarrow_4}  57$  &  $20^{(1)}$
        & $47_{\leftarrow_2} 31_{\leftarrow_4} 41_{\leftarrow_2} 27$  &  $1210^{(3)}$ \\
$8$ & $49_{\leftarrow_4} 65_{\leftarrow_2} 43_{\leftarrow_4} 57$  &  $2120^{(1)}$
        & $\bf{53_{\leftarrow_2} 35_{\leftarrow_2} 23_{\leftarrow_2} 15}$  &  $\bf{1110^{(5)}}$ \\
$9$ & $55_{\leftarrow_4}  73_{\leftarrow_4} 97_{\leftarrow_4} 129$  &  $2220^{(1)}$
        & $59_{\leftarrow_2} 39$  &  $10^{(1)}$ \\
$10$ & $\bf{61_{\leftarrow_4}  81}$  &  $\bf{20^{(3)}}$
        & $65_{\leftarrow_2} 43_{\leftarrow_4} 57$  &  $120^{(1)}$ \\
$11$ & $67_{\leftarrow_4}  89_{\leftarrow_2} 59_{\leftarrow_2} 39$  &  $2110^{(1)}$
        & $71_{\leftarrow_2} 47_{\leftarrow_2} 31_{\leftarrow_4} 41_{\leftarrow_2} 27$  &  $11210^{(5)}$ \\
$12$ & $73_{\leftarrow_4}  97_{\leftarrow_4} 129$  &  $220^{(1)}$
        & $\bf{77_{\leftarrow_2} 51}$  &  $\bf{10^{(5)}}$ \\
$13$ & $79_{\leftarrow_4}  105$  &  $20^{(5)}$
        & $83_{\leftarrow_2} 55_{\leftarrow_4} 73_{\leftarrow_4} 97_{\leftarrow_4} 129$  &  $12220^{(1)}$ \\
$14$ & $\bf{85_{\leftarrow_4}  113_{\leftarrow_2} 75}$  &  $\bf{210^{(1)}}$
        & $89_{\leftarrow_2} 59_{\leftarrow_2} 39$  &  $110^{(1)}$ \\
$15$ & $91_{\leftarrow_4}  \ldots $  &  $221111210^{(3)}$
        & $95_{\leftarrow_2} 63$  &  $10^{(3)}$ \\
$16$ & $97_{\leftarrow_4}  129$  &  $20^{(1)}$
        & $\bf{101_{\leftarrow_2} 67 _{\leftarrow_4} \cdots }$  &  $\bf{12110^{(3)}}$ \\
\hline
\end{tabular}
\end{center}
\end{table}

\subsection{Prefix Properties} \label{subsec:prefix}
Yet another way of analyzing sequence cyclic behavior is to characterize primitive sequences as subsequences of prefixes. Consider Table  \ref{table:prefixes} where it can be seen that there are many instances where one sequence is a prefix of another sequence. In this section, we present a set of rules for generating sets of all prefixes in an ordered fashion. That this is a possibility should be clear from the cyclic nature of the sequences generated as a function of the index $p$ in Table  \ref{table:prefixes}. Appendices A and B form the foundation and provide the details details behind the prefix rules presented in this section.

\begin{observation} \label{obs:wvlttree}
\begin{enumerate}
The following steps can generate the set of all possible prefixes for all possible primitive sequences:
\item (Initial condition) Begin with the sequence $m=0,1,2,3, \ldots$ that imply a sequence of 0-nodes of the form $6m+3$
\item Mark all values of $m$ whose $6m+3$ value will map to $4k+1$ pumps (these values are periodic at every fourth value of $m$).
\item Create new tree levels  as follows.
\begin{itemize}
\item[i.] Create a triple by sequentially taking three values at a time (do not include the marked values as they represent those prefixes that terminate
at a $4k+1$ pump - the end of a primitive sequence).
\item [ii.] Create next prefix ending in a 1-node by taking every other value in the sequence of triples
\item [iii.] Create next prefix ending in a 2-node by taking every middle value in the sequence of triples
\end{itemize}
\item return to step 2
\end{enumerate}
\end{observation}

The above procedure can be made clear via the following demonstration (the numbers shown correspond to values of $m$ as developed in Appendices A and B).
Consider the first level to generate all possible prefixes of length $2$:
\begin{equation}
\begin{matrix}
 & & 01: \underline{0} \; 2 \; 4\; 6\; \underline{8}\; 10 \; \ldots\\
 &\nearrow_1 & \\
0-node: 0 \; 1 \; 2\; \underline{3} \; 4\; 5 \; 6 \; \underline{7} \; 8 \; \ldots  & & \\
 &\searrow_2 & \\
& & 02:1 \; 5 \; 9 \; \underline{13}\; 17 \; 21 \; \ldots\\
\end{matrix}
\end{equation}
Iterate the above procedure to generate all possible prefixes of length $3$:
\begin{equation}
\begin{matrix}
& \nearrow_1  & 011: 2 \; 6\;  10 \; \underline{14} \; 18 \; 22 \; 26 \; \underline{30} \; \ldots \\
01: \underline{0} \; 2 \; 4\; 6\; \underline{8}\; 10 \; \ldots  & \\
 &\searrow_2  & 012: 4 \; \underline{12} \; 20 \; 28 \; \ldots\\
 & & \\
& \nearrow_1  & 021: 1 \; 9 \; 17 \; \underline{25}  \; 33 \; \ldots \\
02:1 \; 5 \; 9 \; \underline{13}\; 17 \; 21 \; \ldots & \\
 &\searrow_2  & 022: 5 \; 21 \;37 \; \underline{53} \; \ldots \\
 \end{matrix}
\end{equation}
Be careful to synchronize this procedure with respect to the $4k+1$ pump. There are various possibilities for iterating this process as the relative
position of sequences leading to $4k+1$ pumps will differ for each tree node. The general form of the generated sequences can be of four possible forms:
\begin{equation*}
\begin{matrix}
1 & 2 & 1 & \underline{X} &  1 & 2 & 1 & \underline{X} & \ldots \\
2 & 1 & \underline{X} &  1 & 2 & 1 & \underline{X} & 1 &\ldots \\
1 & \underline{X} &  1 & 2 & 1 & \underline{X} & 1 & 2 & \ldots \\
\underline{X} &  1 & 2 & 1 & \underline{X} & 1 & 2 & 1 & \ldots \\

\end{matrix}
\end{equation*}

The above procedure can be iterated to map every 0-node to its $4k+1$ pump along a primitive sequence. Observe, the $m$ value is enough
to generate this information. The astute reader will recognize that the above decimation process is analogous to a wavelet transformation thus, once again, hinting at a renormalization procedure. Additionally, this symbolic formulation bears some resemblance to encoding vectors developed for characterizing the convergent behavior of the Collatz function \cite{terras1976}.  However, this work is focused primarily on the quasi-global behavior of families of primitive itineraries having the same symbolic sequence. Along with the C-ladder, primitive itineraries are crucial to understand because they begin at a 0-node, visit only hubs that either $1$-nodes or $2$-nodes and terminate at a $4k+1$ pump of a C-ladder. 

\subsubsection{An Expansive Sequence Example} \label{sec:revisitstability}
While there are many results to consider in light of Section \ref{subsec:symostability}, we once again focus on applying this formulation to 
the most expansive case: generating a sequence of 1-nodes of arbitrary length. Such a construction amounts to traversing the tree outlined in Observation \ref{obs:wvlttree} along its outermost branches (corresponding to a sequence of 1-node transitions). In this case, it can shown that the position of $4k+1$ pump end points at each tree level will alternate between the $0^{th}$ position and the $3^{rd}$ position (e.g. in the above tree figures, consider the position of $4k+1$ pumps for the $10$ case versus the $110$ case).

Given this construction, it is straightforward to see that there always exists an $m$, where, starting at a 0-node value $6m+3$,  it is possible to generate a sequence of 1-nodes of arbitrary length.  However, given the alternating nature of the $4k+1$ end points, once the
expansive sequence encounters the $4k+1$ pump value, the question arises as to whether or not the entry point into the next primitive sequence will be enough to cancel the expansion. We take up this question in part II of this work which addresses the global behavior of the Collatz function.

\section{Conclusions}
In this work, we have presented a symbolic dynamical formalism for characterizing all itineraries generated by the Collatz function. This symbolic approach enables the establishment of dynamical similarities within the presented tree structure. The demonstrated techniques allow for understanding the local and quasi-global behavior of all trajectories. Future work will involve applying this approach to characterize the global behavior of itineraries.


\section*{Disclosure statement}
No potential conflict of interest was reported by the author.

\appendix

\section{Prefix Analysis} \label{sec:AppendixA}

\subsection{Prefixes of length 2}
\subsubsection{The $10$ prefix}

Start with a 0-node $x=6m+3$ and applying
\begin{equation*}
\frac{3x+1}{2}=\frac{3(6m+3)+1}{2}=9m+5
\end{equation*}
to the constraint that the output must be a 1-node
\begin{equation*}
9m+5 = 6p+5
\end{equation*}
leads to the condition
\begin{equation}
p = \frac{3}{2}m.
\end{equation}
This generates the following progression
\begin{equation*}
\begin{matrix}
k = 0,1,2,3,4, \ldots & m=2k = 0,2,4,6,8, \ldots & 6m+3 = 3,15,27,39, 51 \ldots \\
& p = 3k = 0, 3, 6, 9, 12 \ldots  & 6p+5 = \underline{5}, 23, 41, 59, \underline{77} \ldots \\
\end{matrix}
\end{equation*}
where each $6m+3$ 0-node maps to its corresponding $6p+5$ 1-node (4k+1 pumps are underlined and
mark the end point of a primitive sequence). The reader can check that these results are
consistent with Table  \ref{table:prefixes}.

\subsubsection{The $20$ prefix}
Start with a 0-node $x=6m+3$ and applying
\begin{equation*}
\frac{3x+1}{4}=\frac{3(6m+3)+1}{4}= \frac{18m+10}{4}
\end{equation*}
to the constraint that the output must be a 2-node
\begin{equation*}
\frac{18m+10}{4} = 6p+1
\end{equation*}
leads to the condition
\begin{equation}
p = \frac{3}{4}m + \frac{1}{4}.
\end{equation}
This generates the following progression
\begin{equation*}
\begin{matrix}
k = 0,1,2,3,4, \ldots & m=4k+1 = 1,5,9,13,17 \ldots & 6m+3 = 9, 33, 57, 81, 105 \ldots \\
& p = 3k+1 = 1, 4, 7, 10, 13, \ldots  & 6p+1 = 7, 25, 43, \underline{61}, 79 \ldots \\
\end{matrix}
\end{equation*}
where each $6m+3$ 0-node maps to its corresponding $6p+1$ 2-node (4k+1 pumps are underlined and
mark the end point of a primitive sequence). The reader can check that these results are
consistent with Table  \ref{table:prefixes}.

\subsection{Prefixes of length 3}
\subsubsection{The $110$ prefix}

Start with a 0-node $x=6m+3$ and applying
\begin{equation*}
\frac{3 [ \frac{3x+1}{2} ] +1}{2} = \frac{3 [ 9m+5 ] +1}{2} = \frac{27m+16}{2}
\end{equation*}
to the constraint that the output must be a 1-node
\begin{equation*}
 \frac{27m+16}{2} = 6p+5
\end{equation*}
leads to the condition (it will shortly become clear why we leave the fraction
in non-reduced form)
\begin{equation}
p = \frac{9}{4}m + \frac{2}{4}.
\end{equation}
This generates the following progression
\begin{equation*}
\begin{matrix}
k = 0,1,2,3,4, \ldots & m=4k+2 = 2,6,10,14,18, \ldots & 6m+3 = 15, 39, 63, 87, 111, \ldots \\
& p = 9k+ 5  = 5, 14, 23, 32, 41, \ldots  & 6p+5 = 35, 89, 143, \underline{197}, 251, \ldots \\
\end{matrix}
\end{equation*}
where each $6m+3$ 0-node maps to its corresponding $6p+5$ 1-node (4k+1 pumps are underlined and
mark the end point of a primitive sequence). The reader can check that these results are
consistent with Table  \ref{table:prefixes}.

\subsubsection{The $210$ prefix}
Start with a 0-node $x=6m+3$ and applying
\begin{equation*}
\frac{3 [ \frac{3x+1}{2} ] +1}{4} = \frac{ 27m+16}{4}
\end{equation*}
to the constraint that the output must be a 2-node
\begin{equation*}
\frac{ 27m+16}{4}   = 6p+1
\end{equation*}
leads to the condition (it will shortly become clear why we leave the fraction
in non-reduced form)
\begin{equation}
p = \frac{9}{8}m + \frac{4}{8}.
\end{equation}
This generates the following progression
\begin{equation*}
\begin{matrix}
k = 0,1,2,3,4, \ldots & m=8k+ 4 = 4, 12, 20, 28, \ldots & 6m+3 = 27, 75, 123, 171, \ldots \\
& p = 9k+ 5  = 5, 14, 23, 32,  \ldots  & 6p+1 = 31, \underline{85}, 139, 187, \ldots \\
\end{matrix}
\end{equation*}
where each $6m+3$ 0-node maps to its corresponding $6p+1$ 2-node (4k+1 pumps are underlined and
mark the end point of a primitive sequence).

\subsubsection{The $120$ prefix}
Start with a 0-node $x=6m+3$ and applying
\begin{equation*}
\frac{3 [ \frac{3x+1}{4} ] +1}{2} = \frac{54m+34}{8}
\end{equation*}
to the constraint that the output must be a 1-node
\begin{equation*}
\frac{54m+34}{8} = 6p+5
\end{equation*}
leads to the condition
\begin{equation}
p = \frac{9}{8}m - \frac{1}{8}.
\end{equation}
This generates the following progression
\begin{equation*}
\begin{matrix}
k = 0,1,2,3,4, \ldots & m= 8k+1 = 1,9,17,25, \ldots & 6m+3 = 9, 57, 105, 153,  \ldots \\
& p = 9k + 1  = 1,10,19,28, \ldots  & 6p+5 = 11, 65, 119, \underline{173}, \ldots \\
\end{matrix}
\end{equation*}
where each $6m+3$ 0-node maps to its corresponding $6p+5$ 1-node (4k+1 pumps are underlined and
mark the end point of a primitive sequence).

\subsubsection{The $220$ prefix}
Start with a 0-node $x=6m+3$ and applying
\begin{equation*}
\frac{3 [ \frac{3x+1}{4} ] +1}{4} = \frac{54m+34}{16}
\end{equation*}
to the constraint that the output must be a 2-node
\begin{equation*}
\frac{ 27m+17}{8}   = 6p+1
\end{equation*}
leads to the condition
\begin{equation}
p = \frac{9}{16}m + \frac{3}{16}.
\end{equation}
This generates the following progression
\begin{equation*}
\begin{matrix}
k = 0,1,2,3,4, \ldots & m=16k+ 5 = 5, 21, 37, \ldots & 6m+3 = 33, 129, 225, \ldots \\
& p = 9k+ 3  = 3, 12, 21,  \ldots  & 6p+1 = 19, \underline{73}, 127, \ldots \\
\end{matrix}
\end{equation*}
where each $6m+3$ 0-node maps to its corresponding $6p+1$ 2-node (4k+1 pumps are underlined and
mark the end point of a primitive sequence).

\subsection{Prefixes of length 4}
A similar set of basic algebraic steps can be used to generate formulae for
higher order sequences. After a certain point, this exercise becomes unnecessary as a
reproducible pattern arises that can be used to unveil the cyclic nature of the symbolic sequences.

\subsubsection{The $1110$ prefix}
\begin{equation}
p = \frac{27}{8}m + \frac{10}{8}.
\end{equation}
This generates the following progression
\begin{equation*}
\begin{matrix}
k = 0,1,2,3,4, \ldots & m=8k+ 2 = 2,10,18, \ldots & 6m+3 = 15, 63, 111, \ldots \\
& p = 27k+ 8  = 8, 35, 62,  \ldots  & 6p+5 = \underline{53}, 215, 377, \ldots \\
\end{matrix}
\end{equation*}

\subsubsection{The $2110$ prefix}
\begin{equation}
p = \frac{27}{16}m + \frac{14}{16}.
\end{equation}
This generates the following progression
\begin{equation*}
\begin{matrix}
k = 0,1,2,3,4, \ldots & m=16k+ 6 = 6,22,38, \ldots & 6m+3 = 39, 135, 231, \ldots \\
& p = 27k+ 11  = 11,38,65,  \ldots  & 6p+1 = 67, \underline{229}, 391, \ldots \\
\end{matrix}
\end{equation*}

\subsubsection{The $1210$ prefix}
\begin{equation}
p = \frac{27}{16}m + \frac{4}{16}.
\end{equation}
This generates the following progression
\begin{equation*}
\begin{matrix}
k = 0,1,2,3,4, \ldots & m=16k+ 4 = 4,20,36, \ldots & 6m+3 = 27, 123, 219, \ldots \\
& p = 27k+ 7  = 7,34,61,  \ldots  & 6p+5 = 47, 209,371, \ldots \\
\end{matrix}
\end{equation*}

\subsubsection{The $2210$ prefix}
\begin{equation}
p = \frac{27}{32}m + \frac{12}{32}.
\end{equation}
This generates the following progression
\begin{equation*}
\begin{matrix}
k = 0,1,2,3,4, \ldots & m=32k+ 28 = 28,60,92, \ldots & 6m+3 = 171, 363, 555, \ldots \\
& p = 27k+ 24  = 24,51,78,  \ldots  & 6p+1 = 145, 307, 469, \ldots \\
\end{matrix}
\end{equation*}

\subsubsection{The $1120$ prefix}
\begin{equation}
p = \frac{27}{16}m + \frac{5}{16}.
\end{equation}
This generates the following progression
\begin{equation*}
\begin{matrix}
k = 0,1,2,3,4, \ldots & m=16k+ 1 = 1,17,33, \ldots & 6m+3 = 9,105,291, \ldots \\
& p = 27k+ 2  = 2,29,56,  \ldots  & 6p+5 = 17,179,341, \ldots \\
\end{matrix}
\end{equation*}

\subsubsection{The $2120$ prefix}
\begin{equation}
p = \frac{27}{32}m + \frac{13}{32}.
\end{equation}
This generates the following progression
\begin{equation*}
\begin{matrix}
k = 0,1,2,3,4, \ldots & m=32k+ 9 = 9,41,73, \ldots & 6m+3 = 57,249,441, \ldots \\
& p = 27k+ 8  = 8,35,62,  \ldots  & 6p+1 = 49,211,\underline{373}, \ldots \\
\end{matrix}
\end{equation*}

\subsubsection{The $1220$ prefix}
\begin{equation}
p = \frac{27}{32}m - \frac{7}{32}.
\end{equation}
This generates the following progression
\begin{equation*}
\begin{matrix}
k = 0,1,2,3,4, \ldots & m=32k+ 5 = 5,37,69, \ldots & 6m+3 = 33,225,417, \ldots \\
& p = 27k+ 4  = 4,31,58,  \ldots  & 6p+5 = \underline{29},191,353, \ldots \\
\end{matrix}
\end{equation*}

\subsubsection{The $2220$ prefix}
\begin{equation}
p = \frac{27}{64}m + \frac{9}{64}.
\end{equation}
This generates the following progression
\begin{equation*}
\begin{matrix}
k = 0,1,2,3,4, \ldots & m=64k+ 21 = 21,85,149, \ldots & 6m+3 = 129,513,897, \ldots \\
& p = 27k+ 9  = 9,36,63,  \ldots  & 6p+1 = 55, 217, 379, \ldots \\
\end{matrix}
\end{equation*}

\section{Prefix Equations} \label{sec:AppendixB}
Once a sequence begins from a 0-node, it can thereafter be viewed as a binary sequence
consisting of the symbols$\{1,2\}$. The results of Appendix A can therefore be phrased
in terms of binary tree. The head node on the binary tree is a 0-node which can transition to one of two possibilities
\begin{equation}
\begin{matrix}
 & & 01: p= \frac{3}{2}m \\
 &\nearrow_1 & \\
0-node & & \\
 &\searrow_2 & \\
& & 02: p= \frac{3}{4}m + \frac{1}{4}\\
\end{matrix}
\end{equation}
which is the initial condition for the tree.

Transitioning from a given node to a 1-node necessarily introduces a scale factor of $3/2$ which must be applied to the $p$ equation.
Transitioning from a given node to a 2-node necessarily introduces a scale factor of $3/4$ which must be applied to the $p$ equation.
After applying the appropriate scale factor, for any node, a shift value must also be added in order to obtain the
$p$ equation for the next node. After the initial condition, there are four possible shift values from one tree node to the next.
\begin{equation}
\begin{matrix}
1 \rightarrow 1 & +\frac{1}{2} \\
1 \rightarrow 2 & +\frac{1}{2}  \\
2 \rightarrow 1 & -\frac{1}{2}  \\
2 \rightarrow 2 & 0 \\
\end{matrix}
\end{equation}
For example, to build the next level in the tree representing prefixes of length $3$, apply these scaling and shift rules to obtain
\begin{equation}
\begin{matrix}
  & &   &  & 011: p= \frac{9}{4}m + \frac{2}{4} \\
  & &   & \nearrow_1 (+\frac{2}{4}) & \\
  & & 01: p= \frac{3}{2}m  & & \\
  & &   & \searrow_2 (+\frac{4}{8}) & \\
  & &   &  & 012: p= \frac{9}{8}m + \frac{4}{8} \\
 &\nearrow_1 & & & \\
 & & & & \\
0-node & & & & \\
 & & & & \\
  &\searrow_2 & & & \\
  & &   &  & 021: p= \frac{9}{8}m - \frac{1}{8} \\
 & &   & \nearrow_1 (-\frac{2}{4}) & \\
& & 02: p= \frac{3}{4}m + \frac{1}{4} & & \\
 & &   & \searrow_2 (+0) & \\
  & &   &  & 022: p= \frac{9}{16}m + \frac{3}{16} \\
\end{matrix}
\end{equation}
where the length $3$ prefix equations are recovered from the tree.
This process can be iterated to recover any $p$ equation for any
primitive sequence emanating from a $0-node$.

As pointed out in the main text of Section \ref{subsec:prefix}, there are important questions
regarding divergent sequences that can be answered without the need for the complete
set of $p$ equations for all possible binary sequences.  One important point to consider for the
$p$ values would be to ensure that the complete set would recover all possible odd integers $\bigcup_p \{6p+1,6p+5\}$.
If it can be shown that every nonnegative integer can be generated by the complete set of $p$ values,
then it can be concluded that every 1-node and 2-node exists in the path of some primitive itinerary.

\end{document}